\documentclass[titlepage,12pt]{article}
\usepackage{amssymb}
\input psfig.sty

\newcommand{\mle}{{MLE}}
\newcommand{\mele}{{MELE}}
\newcommand{\mse}{{MSE}}
\newcommand{\PMC}{{\rm PMC}}

\title{\sc Mean Likelihood Estimators}
\author{A.I. McLeod and B. Quenneville
\thanks{A.I. McLeod is Professor,
Department of Statistical and Actuarial Sciences,
University of Western Ontario, London, Ontario, Canada, N6A 5B7,
e-mail: aim@uwo.ca.
B. Quenneville is Senior Methodologist,
Time Series Research and Analysis Centre,
Statistics Canada, Ottawa, Ontario, Canada, K1A~0T6.
e-mail: quenne@statcan.ca.
The authors wish to thank Jamie Stafford and two anonymous referees for their
very helpful and insightful remarks.}}
\date{March 1, 1999}
\begin{document}
\maketitle
\pagenumbering{roman}

\section*{Abstract}
The use of {\it Mathematica} in deriving mean likelihood estimators is
discussed.
Comparisons between the maximum likelihood estimator, the mean
likelihood estimator and the Bayes estimate based on
a Jeffrey's noninformative prior using the criteria
mean-square error and Pitman measure of closeness.
Based on these criteria we
find that for the first-order moving-average time series model,
the mean likelihood estimator outperforms the maximum likelihood
estimator and the Bayes estimator with a Jeffrey's noninformative prior.

{\it Mathematica} was used for symbolic and numeric
computations as well as for the graphical display of results.
A {\it Mathematica} notebook is available which provides
supplementary derivations and code from
http://\-www.stats.uwo.ca/\-mcleod/\-epubs/\-mele.
The interested reader can easily reproduce or extend any
of the results in this paper using this supplement.

{\bf Keywords:}
Binomial Distribution;
Efficient Likelihood Computation;
Exponential Distribution;
First-order moving-average time series model;
Mean Square Error Criterion;
Pitman Measure of Closeness
\clearpage
\pagenumbering{arabic}
\pagestyle{headings}

\section{Introduction}
The maximum likelihood estimator (\mle) is perhaps the most common and widely
accepted estimator of a parameter in a statistical model denoted by
$({\mathcal S},\Omega, f)$, where ${\mathcal S},\Omega, f$ denote respectively
the sample space, the parameter space and the probability density function
(pdf)\@.
We will take ${\mathcal S}= {\mathcal R}^n,  X =
\left(X_1, X_2, \ldots, X_n\right)  \in {\mathcal S}$,
and $f(x,\theta)$\@.
In the standard case of independent and identically distributed observations,
$f(x,\theta)= \Pi_{i = 1}^n f_1(x_i)$,
where $f_1(x)$ is the pdf of $X_1$\@.
Given data $X$, the likelihood function is
$L(\dot{\theta}) = f(X; \dot{\theta}),\; \dot{\theta} \in \Omega$
and the \mle\ of the parameter $\theta$ is defined as that value
$\dot{\theta}$ which globally maximizes $L(\dot{\theta})$\@.
{\it Mathematica} (Wolfram, 1996) has been widely used in the study of
fundamental and general
aspects of maximum likelihood estimation --- see Andrews and Stafford (1993);
Stafford and Andrews (1993);   Stafford,  Andrews and Wang  (1994).
As well {\it Mathematica} has been used for obtaining symbolically exact maximum
likelihood estimators in situations where the use of numerical techniques are
less convenient such as with grouped or censored data or logistic regression ---
see Cabrera  (1989);  Currie (1995).

For simplicity we will deal with the case where $\Omega$ is one-dimensional.
The multidimensional case may in general be reduced to the one-dimensional case
by using marginal, conditional or concentrated likelihoods or by integrating
over the nuisance parameters whichever is more suitable in a particular
situation.  Under the usual regularity conditions, the \mle, $\hat{\theta}$, is
approximately normally distributed with mean $\theta$ and covariance matrix
$I_{\theta}^{-1}$, where  $I_{\theta}$ denotes the Fisher information
matrix.
It is also true that the mean likelihood estimator (\mele) is equally
efficient in large samples.
In general the \mele, $\bar{\theta}$ is defined by
\[
\bar{\theta} =
\frac{\int_{\Omega} \dot{\theta} L(\dot{\theta}) d\dot{\theta}}{\int_{\Omega}
L(\dot{\theta}) d\dot{\theta}},
\]
where $L(\dot{\theta})$ is the likelihood function.
It should be noted that although the \mele\ is identical to the Bayes
estimator with a uniform prior, it is not often considered in
frequentist settings even though
Pitman (1938) showed that when the problem is location invariant, the
\mele\ is the best invariant estimator.
Barnard, Jenkins and Winsten (1962)
recommended the \mele\ for time series problems and suggested that it will
often have lower \mse\ than the \mle.
In changepoint analysis, where the usual regularity conditions for the
\mle\ do not hold and the \mle\ is inefficient but the \mele\ works well
(Ritov, 1990; Rubin and Song, 1995).

Unlike the \mle\ the \mele\ is not invariant under reparameterization.
Although the \mele\ has a Bayesian interpretation, it is not
the Bayesian estimator that is usually recommended.
In order that the estimator share \mle\ property of being invariant
under parameter transformation, the Jeffrey's noninformative
prior is recommended
when there is no prior information available (Box and Tiao, 1973, \S 1.3).
The Jeffrey's prior is given by
$p(\theta) \propto \sqrt I_{\theta}$.

There are situations, such as in
the first-order moving-average model (MA(1)) where the \mle\ in
finite samples has non-zero probability of lying on the boundary of the
parameter region but this phenomenon does not happen with the \mele\
or Bayesian estimator as can be seen from the following result.

{\bf Theorem 1:} Let $\Omega = \left[a,b\right]$ then
$\Pr \left\{ \bar{\theta} \in (a,b) \right\}  = 1$\@.

{\bf Proof:}
The likelihood function, $L(\dot{\theta})$,  defined below, is easily seen to
be continuous and differentiable in the interval $[a, b]$ and non-negative.
It then follows from the generalized mean-value theorem
(Borowski and Borwein, 1991, p.371) that
$\bar{\theta} \in (a,b).\;    \square$

In many cases the \mle\ is easy to compute using pen and paper.
However with {\it Mathematica} we can now easily obtain the \mele\ by numerical
integration and sometimes symbolically.
In fact, for problems where the likelihood function is complicated or difficult
to evaluate the \mele\ may be computationally easier to compute than the
traditional \mle. As shown in Theorem 2, both the \mle\ and \mele\ are first order
efficient.

{\bf Theorem 2:}
Under the usual regularity conditions for maximum likelihood estimators,
$\bar{\theta}  = \hat{\theta} + O_p(1/n)$\@.

{\bf Proof:}
The likelihood function, $L(\dot{\theta})$, is to $O_p(1/n)$ equal to the
normal probability density function with mean $\theta$ and variance
$I_{\theta}^{-1}$ (Tanner, 1993, p.16).
The result then follows directly from this approximation.   $\square$

Now consider an estimator $\hat{\theta_1}$ of $\theta$.
The mean-square error (\mse) of an estimator $\hat{\theta_1}$  is defined as
$\sigma^2(\hat{\theta_1} | \theta ) =  E \left\{ ( \hat{\theta_1} - \theta )^2
\right\}$\@.
The relative efficiency of $\hat{\theta_1}$ vs $\hat{\theta}$ is defined as
$R(\hat{\theta_1}, \hat{\theta} | \theta )  = \sigma^2(\hat{\theta} | \theta ) /
\sigma^2(\hat{\theta_1} | \theta )$\@.
Clearly, from Theorem 2, as
$n\rightarrow  \infty,\; R(\bar{\theta}, \hat{\theta} | \theta ) = 1$\@.
Barnard, Jenkins and Winsten (1962) suggested that in many small sample
situations the \mele\ is preferred by the mean-square error criterion
and hence at least for some values of $\theta$,
$R(\bar{\theta}, \hat{\theta} | \theta ) > 1$,
where $\hat{\theta}$ and $\bar{\theta}$ denote  the \mle\ and \mele\ respectively.

Pitman (1937) formulated a useful alternative to the \mse\ in the situation
where no explicit loss function is known.
Consider two estimators,
$\hat{\theta_1}$ and $\hat{\theta_2}$,
and assume that with probability one,
$\hat{\theta_1} \ne \hat{\theta_2}$ then
the Pitman measure of closeness
for comparing $\hat{\theta_1}$ vs $\hat{\theta_2}$ is defined as
\begin{equation}
{\rm PMC} \left[ \hat{\theta_1}, \hat{\theta_2} | \theta \right] =
\Pr \left\{ | \hat{\theta_1}-\theta | < | \hat{\theta_2}-\theta | \right\}.
\label{PMCone}
\end{equation}
When
${\rm PMC} \left[ \hat{\theta_1}, \hat{\theta_2} | \theta \right] > 1/2$,
$\hat{\theta_1}$ is preferred to $\hat{\theta_2}$.
The monograph of Keating, Mason and Sen (1993) provides an extensive
survey of recent work and applications of the PMC.
Additionally, volume 20 (11) of
{\it Communications in Statistics: Theory and Methods\/}
contains an entire issue on the PMC.

Unlike the \mse\ and relative efficiency, the PMC depends on the bivariate
distribution of the two estimators.
The PMC is more appropriate in many scientific and industrial applications
in which the estimator which is closer to the truth is required.
Sometimes it is felt that the \mse\ and other risk criteria give too much weight
to large deviations which may seldom occur.
Rao and other researchers (Keating, Mason and Sen, 1993, \S3.3) have found that
risk functions such as \mse\ and mean-absolute-error can often be shrunk but that
this shrinkage occurs at the expense of the PMC.
The \mse\ or some other risk function is more appropriate than PMC in the
decision theory framework when there is some economic or other loss associated
with the estimation error.
In practice it is often useful to consider both the PMC and \mse\ and in many
situations there appears to be a high level of concordance between these
estimators (Keating, Mason and Sen, 1993, \S2.5).

As originally pointed by Pitman (1937) the PMC criterion is intransitive
but it is arguable whether this is a practical limitation.
This point as well as other limitations and extensions of the PMC are
discussed by Keating, Mason and Sen (1993, Ch.3)

{\bf Theorem 3:}
$\bar{\theta}$ and $\hat{\theta}$
are not necessarily asymptotically equivalent under the PMC.

{\bf Proof:}
See eqn. \ref{eqBinomialLimit}.  $\square$

The next theorem shows that the \mele\ minimizes the mean likelihood of the
squared error.

{\bf Theorem 4:}
Choosing $\dot{\theta}  = \bar{\theta}$ minimizes $\delta(\dot{\theta})$,
where
\[
\delta(\dot{\theta})= \int_\Omega \left( \dot{\theta}-
\theta \right)^2 L \left( \dot{\theta} \right) d\dot{\theta}.
\]

{\bf Proof:}
Using calculus, the result follows directly.  $\square$

{\bf Theorem 5:}
$\bar{\theta}$ is a function of the sufficient statistic for $\theta, S$,
if there is one.

In general, the \mele\ is a biased estimator.

{\bf Theorem 6:}
If $\Omega$ has compact support and $0 < {\rm Var}( \bar{\theta} )< \infty$
then $ E\left\{\bar{\theta}\right\} \not= \theta$.

Theorems 5 and 6 are derived in Quenneville (1993).
The \mele\ is formally equivalent to the Bayes estimator under a locally uniform
prior with the squared error risk function and many of the above theorems have
their well-known Bayesian analogues.

We are now going to make comparisons between these three estimators
for three statistical models: Bernouilli trials, exponential lifetimes
and the first-order moving average process.
The symbolic, numeric and graphical computations will all be done using
{\it Mathematica\/}.
The interested reader can reproduce or extend our computations using
the {\it Mathematica\/} notebooks we have provided
(McLeod and Quenneville, 1999).
Frequentist analysis of Bayesian estimators is not often done but
Dempster (1998) and Quenneville and Singh (1999)
have argued that frequentist considerations are obviously informative
even in the Bayesian setting.

\section{Bernoulli Trials}

\begin{figure}
\centerline{\psfig{figure=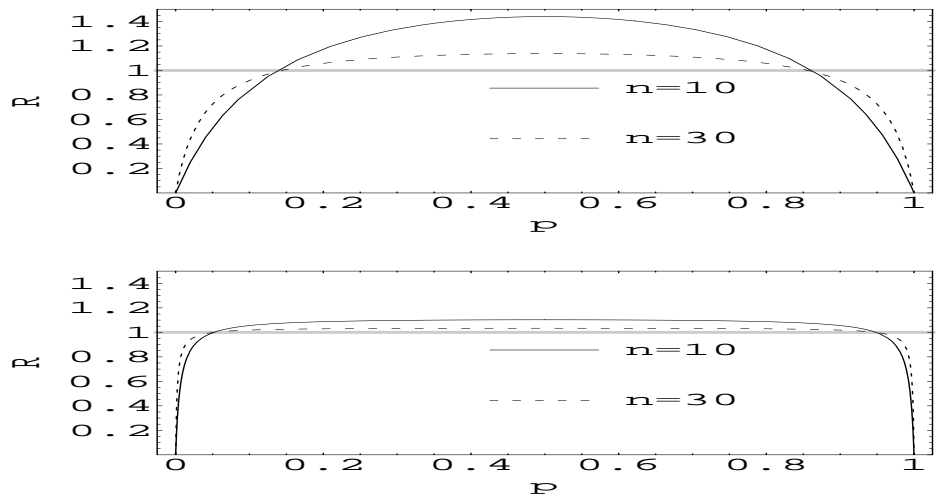,width=5.5in,height=2.75in}}
\caption{
Relative efficiency of alternative binomial estimators.
Top panel: \mele, relative efficiency, $R(\bar{p}, \hat{p} | p)$ for $n=10, 30$.
Bottom panel: Bayes estimator with Jeffrey's prior,
relative efficiency,
$R(\tilde{p}, \hat{p} | p)$ for $n=10, 30$.
}
\label{BinomialMSE}
\end{figure}

We will now examine the performance of these three estimators in the estimation
of the parameter $p$ in a sequence of $n$ Bernoulli trials where $X$ is the
observed number of successes and $p$ is the probability of success.
The probability function is
\[
f_x(n,p) = \left( \begin{array}{c} n\\x \end{array} \right) p^x (1-p)^{n-x}.
\]
So if $X$ successes are observed in $n$ trials, the likelihood function may be
written $L(p) = p^X(1 - p)^{(n - X)}$ and the \mle\ may be derived by calculus,
$\hat{p} = X/n$\@.
Using {\it Mathematica} it is easily shown that the \mele\ of
$p$ is $\bar{p} = (X + 1)/(n + 2)$ and that $R(\bar{p}, \hat{p} | p) > 1 $
provided
\[
p \in \left(\frac{2n -
\sqrt{2 n^2 + 3n +1} +1}{2 (2n+1)},\frac{2n +
  \sqrt{2 n^2 + 3n +1} +1}{2 (2n+1)} \right).
\]
As shown in Figure~\ref{BinomialMSE}, the \mele\ is always more efficient
over most
of the range and the relative efficiency tends to $1$ as
$n \rightarrow \infty$\@.

It is interesting to compare the \mele\ with Bayes estimate under a
Jeffrey's prior.
The Jeffrey's prior for $p$ is (Box and Tiao, p.35),
$\pi(p) = 1/\sqrt{ p (1-p) }$.
Combining with the likelihood we can use {\it Mathematica\/} to show
that the resulting Bayes estimator is
$\tilde{p} =(1+4 X)/(2+4 n)$.
From Figure~\ref{BinomialMSE},
we see that the Bayes estimator with Jeffrey's prior
tends have smaller mean-square error
over an even slightly larger range of $p$ than the \mele but
the gain in efficiency with the mele can be greater.
As with the \mele, the relative efficiency tends to $1$ as
$n \rightarrow \infty$\@.
Once again,
using {\it Mathematica} we can show that
$R(\tilde{p}, \hat{p} | p) > 1 $
provided
\[
p \in \left(
{\frac{1 + 5\,n - {\sqrt{1 + 9\,n + 20\,{n^2}}}}
   {2\,\left( 1 + 5\,n \right) }},
{\frac{1 + 5\,n + {\sqrt{1 + 9\,n + 20\,{n^2}}}}
{2\,\left( 1 + 5\,n \right) }}
\right).
\]

The PMC criterion given in eqn.~\ref{PMCone}
is not applicable in the case of the
binomial since due to the discreteness there can be ties in the
values of the estimators.
Keating, Mason and Sen (1993, \S 3.4.1) and one of the referees
have suggested the following modified version of Pitman's measure
of closeness,
\[
{\rm PMC} \left[ \bar{\theta}, \hat{\theta} | \theta \right] =
\Pr \left\{ | \bar{\theta}-\theta | < | \hat{\theta}-\theta | \right\}
+ {1\over 2}
\Pr \left\{ | \bar{\theta}-\theta | = | \hat{\theta}-\theta | \right\}.
\]
With this modification, PMC is transitive and reflexive.

\begin{figure}
\centerline{\psfig{figure=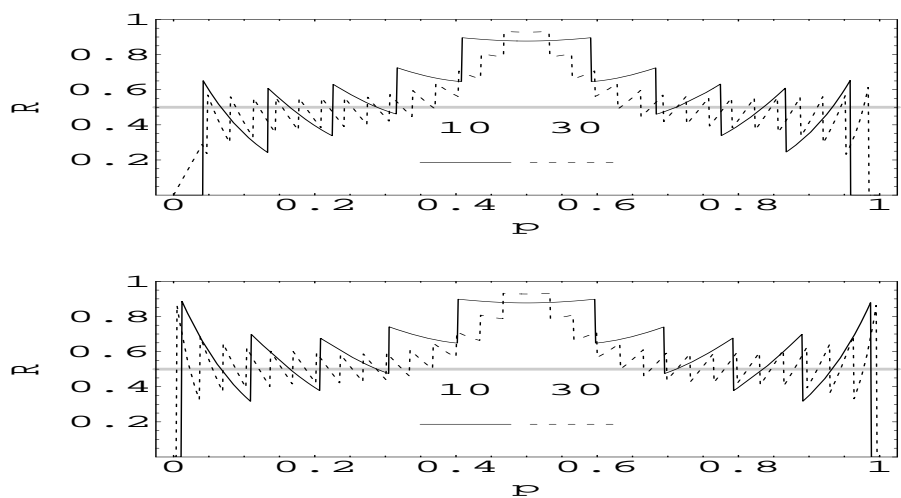,width=5.5in,height=2.75in}}
\caption{
Pitman measure of closeness for alternative binomial estimators.
Top panel: \mele, $\PMC(\bar{p}, \hat{p} | p)$ for $n=10, 30$.
Bottom panel: Bayes estimator with Jeffrey's prior
$\PMC(\tilde{p}, \hat{p} | p)$ for $n=10, 30$.
}
\label{PMCbinomial}
\end{figure}

Figure~\ref{PMCbinomial} suggests
the following asymptotic result,
\begin{equation}
\lim_{n \rightarrow \infty} {\rm PMC}(\bar{p}, \hat{p} | p) =
      \left\{
      \begin{array}{ll}  1           & \mbox{$p=0.5$} \\
                        \frac{1}{2} & \mbox{$p \not= 0.5, 0, 1$} \\
                        0           & \mbox{$p= 0, 1$}
      \end{array}
      \right.
\label{eqBinomialLimit}
\end{equation}
This result may be established using the Geary-Rao Theorem
(Keating, Mason and Sen, p.103).

Figure~\ref{PMCbinomial} also suggests that in terms of the
PMC the advantage over the \mle\ of the \mele\ or of the Bayes estimate with a
Jeffrey's prior disappears when there is no prior information
about $p$.

\section{Exponential Lifetimes}

Consider a sample of size $n$ denoted by $X_1, \ldots ,  X_n$ from an
exponential distribution with mean $\mu$ and let $T = \sum_{i = 1}^{n} X_i$\@.
The likelihood function for $\mu$ can be written
$L(\mu) = \mu^{-n} e^{-T/\mu}$, the \mle\ of $\mu$ is given by
$\hat{\mu} = T/n$ and the \mele\ of $\mu$ is $\bar{\mu} =  T/(n-2)$\@.
The Jeffrey's prior for $\mu$ can be taken as $\mu^{-1}$ which produces
a Bayesian estimate, $\tilde{\mu}=T/(n-1)$.

\begin{figure}
\centerline{\psfig{figure=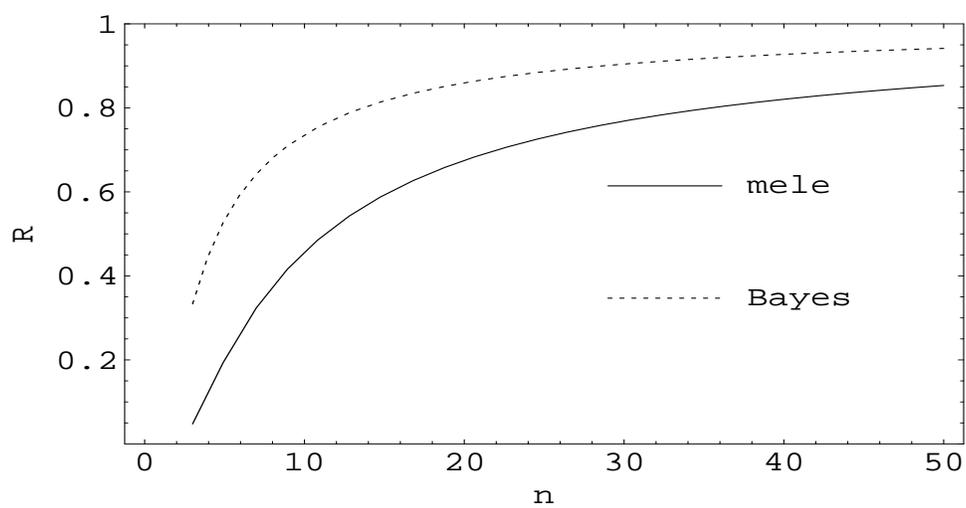,width=5.5in,height=2.75in}}
\caption{
Relative efficiency $R$
of the \mele\ and Bayes estimator vs the \mle\
of the mean $\mu$ in a random sample of
size $n$ from an exponential distribution.
}
\label{Rlife}
\end{figure}

A simple computation with {\it Mathematica\/} gives the
relative efficiency,
\begin{eqnarray*}
R(\bar{\mu}, \hat{\mu}) &=&
{\frac{1}{n}} + {\frac{-5 + n}{4 + n}} \\
 &=& 1 - {\frac{8}{n}} + {\frac{36}{{n^2}}} - {\frac{144}{{n^3}}} +
  {\frac{576}{{n^4}}} - {\frac{2304}{{n^5}}} +
  {{{\mathcal O}({\frac{1}{n}})}^6}.
\end{eqnarray*}
Similarly,
$R(\tilde{\mu}, \hat{\mu})=1+1/n+4/(n+1)$.
Figure~\ref{Rlife} shows that the \mele\ can be much less efficient.

\begin{figure}
\centerline{\psfig{figure=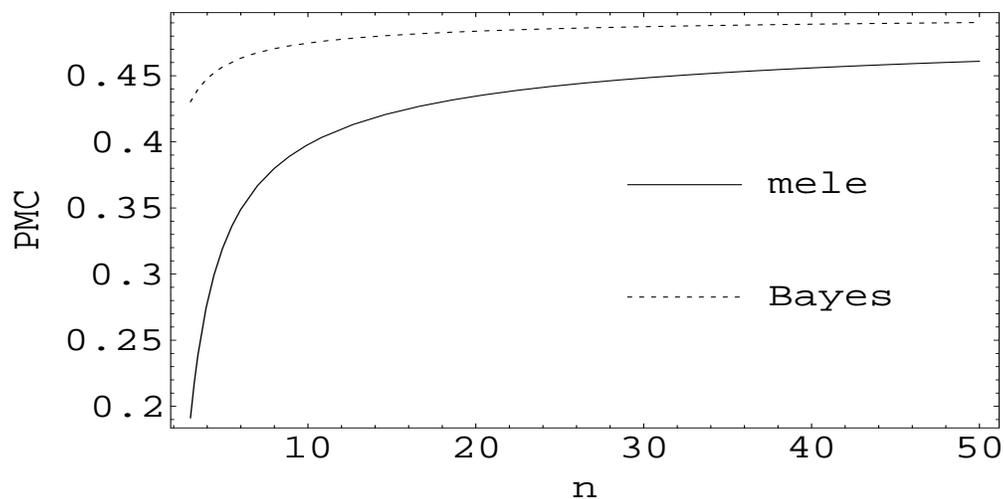,width=5.5in,height=2.75in}}
\caption{
Pitman Measure of Closeness, PMC,
of the \mele\ and Bayes estimator vs the \mle\
of the mean $\mu$ in a random sample of
size $n$ from an exponential distribution.
}
\label{PMClife}
\end{figure}

Since $T$ has a standard gamma distribution with shape parameter $n$ and scale
parameter $\mu$, the PMC is easily evaluated using the Geary-Rao Theorem
(Keating, Mason and Sen, 1993, p.103).
Letting $a=\bar{\mu}$ or $a=\tilde{\mu}$, we can write
\[
{\rm PMC}(a, \hat{\mu} | \mu)=
\int_{0}^{b \mu} \frac{e^{-x/\mu} x^{n-1} \mu^{-n}}{\Gamma(n)} dx
\]
where $b=n(n-2)/(n-1)$ or $b=2 n (n-1)/(2n-1)$ according as
$a=\bar{\mu}$ or $a=\tilde{\mu}$\@.
Notice that without loss of generality we may take $\mu =1$ since
${\rm PMC}(\bar{\mu}, \hat{\mu} | \mu) =
   {\rm PMC}(\bar{\mu}, \hat{\mu} | 1)$\@.
From Figure~\ref{PMClife}, ${\rm PMC}(a, \hat{\mu} | \mu) < 0.5$
for both $a=\bar{\mu}$ or $a=\tilde{\mu}$\@.

It is sometimes mistakenly thought that Theorem 4 or its Bayesian analogue
guarantees that at least over some region of the parameter space, the \mele\
and the Bayes estimator
will have outperform the \mle\
but this need not be the case.

\section{MA(1) Process}
\subsection{Introduction}

The MA(1) time series with mean $\mu$ may be written
$Z_t  = \mu + A_t  + \theta A_{t - 1}$,
where $Z_t$ denotes the observation at time $t=1, 2, \ldots$ and $A_t$,
the innovation at time $t$, is assumed to be a sequence of independent normal
random variables with mean zero and variance $\sigma_A^2$\@.
The parameter $\theta$ determines the autocorrelation structure of the series
and for identifiability we will assume that $|\theta| \le 1$\@.
When $| \theta |< 1$, the model is invertible
(Brockwell and Davis, 1991, \S3.1).
For simplicity we will examine the case where $\mu= 0$\@.
Such MA(1) models often arise in practical applications as the model for a
differenced nonstationary time series.
The non-invertible case $\theta=1$ occurs when a series is over-differenced.

In large-samples, standard asymptotic theory suggests that the maximum
likelihood estimate for $\theta$, denoted by $\hat{\theta}$,  is approximately
normal with mean $\theta$ and variance $(1 - \theta^2)/n$ where $n$ is the
length of the observed time series.
Cryer and Ledolter (1981) established the somewhat surprising result that
$\Pr \{\hat{\theta} = \pm1 \} > 0$\@.
This result holds for all finite $n$ and for all values of $\theta$\@.
For example when $n=50,\; \Pr \{\hat{\theta} = 1 | \theta = 0 \} = 0.002$ and
$\Pr \{\hat{\theta} = 1 | \theta = 0.8 \} = 0.13$
(Cryer and Ledolter, 1981, Table 2).
Let $\bar{\theta}$ denote the mean likelihood estimate of $\theta$\@.
In view of Theorem 1, this problem does not occur with $\bar{\theta}$\@.

Now we will show that the \mele\ dominates the \mle\ both for the \mse\ and PMC
criteria when $n=2$\@.  When $n=50$, the \mele\ is better than the \mle\ unless the
parameter $\theta$ is very close to $\pm 1$\@.
Since even the useless estimator obtained by ignoring the data and setting the
estimate to $1$ does better when $\theta=1$, we can conclude that \mele\ is
generally a better estimator. Further mean-square error computations which
support this conclusion for other values of $n$ are given by Quenneville (1993)
and can be verified by the reader using the electronic supplement.

\subsection{Exact Results for $n=2$}

Given a Gaussian time series of length $2, Z_1, Z_2$, generated from the
first-order moving average equation $Z_t = A_t-\theta A_{t - 1}$,
where $A_t$ are independent normal random variables with mean zero and
variance $\sigma_A^2$\@.  Let $W = -Z_1 Z_2 /(Z_1^2 + Z_2^2)$\@.
Then given data,  $Z_1, Z_2$, the exact concentrated likelihood function for
$\theta$ is (Cryer and Ledolter, 1981; Quenneville, 1993),
\[
L(\theta | W ) =
\frac{\sqrt{1+\theta^2+\theta^4}}{1+\theta^2-2 \theta W}
\]
and
\[
\hat{\theta} = \left\{
\begin{array}{ll}
-1 & \mbox{$W \in [-0.5, -0.25]$}  \\
\frac{1- \sqrt{1-16W^2}}{4W} & \mbox{$W \in (-0.25, 0.25), W \not= 0$} \\
0 & \mbox{$W =0$}  \\
1 & \mbox{$W \in [0.25, 0.5]$.}
\end{array}
\right.
\]

Unfortunately $\bar{\theta}$, cannot be evaluated symbolically.
However using {\bf NIntegrate} we can obtain it numerically.
Numerical evaluation suggests that $\bar{\theta}$ is either a linear or close
to a linear function of $W$.
To speed up our computations for the mean-square error of $\bar{\theta}$,
we use the {\bf FunctionInterpolation} in {\it Mathematica} to construct
$\bar{\theta} = \bar{\theta}(W)$\@.
The \mse\ and PMC for $\bar{\theta}$ and $\hat{\theta}$ are easily evaluated
numerically using the  pdf of $W, \; f_W(x)$, derived by Quenneville (1993),
\[
f_W(x) =
\frac{2 \sqrt{1+\theta^2+ \theta^4}}{\pi \sqrt{1-4x^2}(1+\theta^2 -2 \theta x)},
    \; |x|\le 1/2.
\]

\begin{figure}
\centerline{\psfig{figure=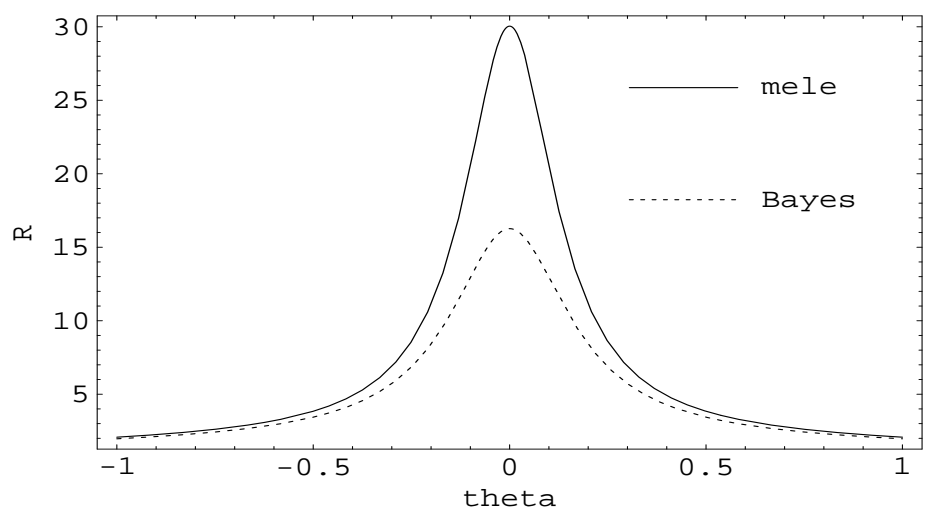,width=5.5in,height=2.75in}}
\caption{
Relative efficiency, $R$,
of \mele\ and Bayes estimator with Jeffrey's noninformative
prior in the MA(1) model with $n=2$.
}
\label{Rma1n2}
\end{figure}

\begin{figure}
\centerline{\psfig{figure=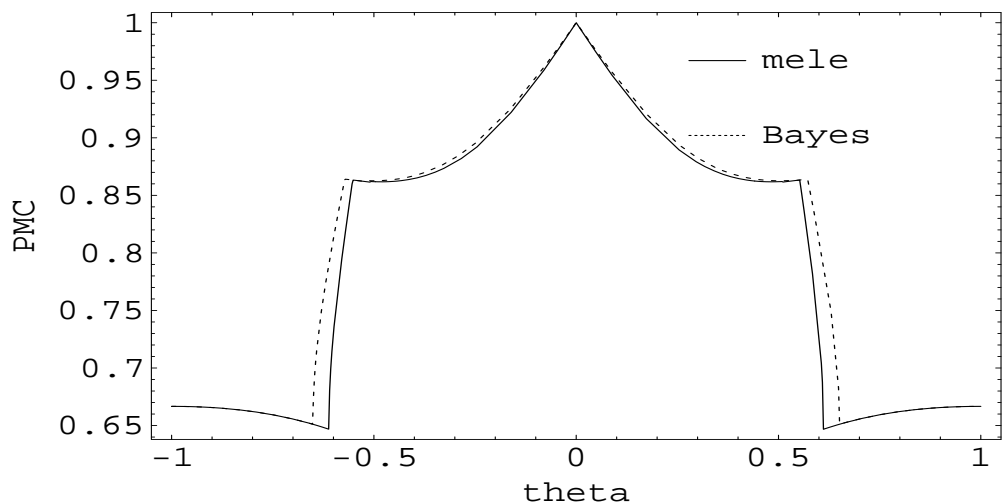,width=5.5in,height=2.75in}}
\caption{
Pitman measure of closeness, PMC,
of \mele\ and Bayes estimator with Jeffrey's noninformative
prior in the MA(1) model with $n=2$.
}
\label{PMCma1n2}
\end{figure}

From Figures~\ref{Rma1n2} and \ref{PMCma1n2},
it is seen that both the \mele\ and Bayesian estimator dominate
the \mle\ both for the \mse\ and PMC criteria.
The \mele\ is slightly better according to the \mse\ but according
to the PMC the Bayes estimator is slightly better than the \mele.

\subsection{Exact Symbolic Likelihood}

Consider the MA(1) process defined by $Z_t  =  A_t  - \theta  A_{t - 1}$,
where $A_t$ is assumed to be normal and independently distributed with mean
zero and variance $\sigma_A^2$\@.
Given $n$ observations $Z^\prime = (Z_1, \ldots, Z_n)$ the exact log likelihood
function of an ARMA process can be written (Newbold, 1974),
\[
\log L(\theta, \sigma_A^2) =
 -\frac{n}{2} \log(\sigma_A^2) - \frac{1}{2} \log(D)-\frac{1}{2\sigma_A^2}
  S(\theta),
\]
where $h^\prime =
(1, \theta, \theta^2, \ldots, \theta^n), \; D= h^\prime h$
and
\[
S(\theta) = (L z - h h^\prime L z/D)^\prime ( Lz - h h^\prime L z/D),
\]
where $L$ is the $(n+1)$ by $n$ matrix,
\[
L = \left(
    \begin{array}{cccccc}
0 & 0 & 0 & \ldots & 0 & 0 \\
1 & 0 & 0 & \ldots & 0 & 0 \\
\theta & 1 & 0 & \ldots & 0 & 0 \\
\theta^2 & \theta & 1 & \ldots & 0 & 0 \\
\vdots & \vdots & \vdots & \ldots & \vdots & \vdots \\
\theta^{n-2} & \theta^{n-3} & \theta^{n-4} & \ldots & 1 & 0 \\
\theta^{n-1} & \theta^{n-2} & \theta^{n-3} & \ldots & \theta & 1
\end{array}
\right).
\]
Maximizing over $\sigma_A^2$ the concentrated log likelihood is given by
\[
\log L_M(\theta) =
 -\frac{n}{2} \log\left[S(\theta)/n \right] - \frac{1}{2} \log(D).
\]
This expression for the concentrated loglikelihood is just as easy to write in
{\it Mathematica} notation as it is in ordinary mathematical notation.
Moreover, it can be evaluated symbolically or numerically.

\begin{verbatim}
LogLikelihoodMA1[t_, z_] :=
    Module[{n = Length[z], Lz, h, detma1, v, Sumsq},
      Lz = Join[{0},
          Table[Sum[z[[i]]] t^(j-i), {i, 1, j}], {j, 1, n}]];
      h = Table[t^j, {j, 0, Length[z]}];
      detma1 = h . h;
      v = -h . Lz/detma1;
      Sumsq = (Lz + h v). (Lz + h v);
      -n/2 Log[Sumsq/n /. t -> t] -
        1/2 Log[detma1 /. t -> t]
    ];
\end{verbatim}

\subsection{Efficient Numeric Likelihood Computations}

Newbold's algorithm can be made much more efficient when only numerical values
of the log likelihood are needed by using the {\it Mathematica} Compiler and
by re-writing the calculations involved to make more use of efficient
{\it Mathematica} functions such as {\bf NestList}, {\bf FoldList} and
{\bf Apply}.
First consider the computation of the vector $L z$ which is of length $n+1$\@.
After some simplifications, we see that $L z  = (\alpha_j)^\prime$, where
$\alpha_0 = 0$ is the first element and the remaining elements are defined
recursively by
$\alpha_j = \theta  \alpha_{j - 1} + Z_j, \;  j = 1, 2, \ldots, n$,
where $Z_0 = 0$\@.
This computation is efficiently performed by {\it Mathematica}'s
{\bf FoldList}\@.
When we are just interested in numerical evaluation we use the compile function
to generate native code which runs much faster.
\begin{verbatim}
GetLz=Compile[{{t,_Real},{z, _Real, 1}},
  FoldList[(#1 t + #2)&,0,z]];
\end{verbatim}

The determinant, $D = 1 + \theta^2 + \theta^4 + \ldots + \theta^{2n}$, is
efficiently computed using {\bf NestList} to generate the individual terms and
then summing.
\begin{verbatim}
DetMA =Compile[{{t,_Real},{n, _Integer}},
  Apply[Plus,NestList[#1 t &,1,n]^2]];
\end{verbatim}

Next, we evaluate the term $h L z/D$\@.
Since
$h L z  = \theta \alpha_1 + \theta^2 \alpha_2 + \ldots + \theta^n \alpha_n$
we can use Horner's Rule to efficiently compute this sum.
Horner's Rule is implemented in {\it Mathematica} using the function {\bf Fold}.
\begin{verbatim}
Getu0 =Compile[{{t,_Real},{Lz, _Real, 1},{detma, _Real}},
       -Fold[#1 t + #2&,0,Reverse[Lz]]/detma];
\end{verbatim}

The computation of the sum of squares function $S(\theta)  = (L z  - h h^\prime L z / D)^\prime (L z  - h h^\prime L z / D)$ is straightforward. The {\it Mathematica} compiler can be used to optimize the vector computations.
\begin{verbatim}
GetSumSq =
  Compile[{{t,_Real},{Lz, _Real, 1},{u, _Real},{n, _Integer}},
    Apply[Plus,(Lz+NestList[#1 t &,1,n] u)^2]];
\end{verbatim}

Finally, the concentrated loglikelihood function is defined.
The computation speed is increased by about a factor of $50$ times when $n=50$
and is even larger for larger $n$\@.
\begin{verbatim}
logLMA1F[t_, z_] :=
  Module[{n=Length[z]},
    Lz=GetLz[t,z];
    detma=DetMA[t,n];
    u=Getu0[t,Lz,detma];
    S=GetSumSq[t,Lz,u,n];
    -(1/2) Log[detma]- (n/2)Log[S/n]];
\end{verbatim}

This function can be maximized using {\it Mathematica\/}'s
nonlinear optimization function {\bf FindMinimum}.

The mean likelihood estimate $\bar{\theta}$ can be evaluated using
{\bf NIntegrate}.
\begin{verbatim}
Meanle[z_]:=
  NIntegrate[t E^logLMA1F[t, z],{t,-1,1}]/
    NIntegrate[E^logLMA1F[t, z],{t,-1,1}]
\end{verbatim}

Notice that in the above expression the loglikelihood function is evaluated
separately in both the numerator and denominator.
Hence, we can save function evaluations by using our own numerical
quadrature routine.
\begin{verbatim}
SimpsonQuadratureWeights[k_,a_, b_]:=
  With[{h=(2 k)/3},
   {a+(b-a)Range[0,2 k]/(2k),
   Prepend[Append[Drop[Flatten[Table[{4,2},{k}]],{-1}],1],1]}]

{X,W}=SimpsonQuadratureWeights[100,-1,1];

GETMEANLEF=
  Compile[{{z, _Real, 1},
           {W, _Real, 1},{X, _Real, 1},{f, _Real, 1}},
           Plus@@(X f)/Plus@@f];			
MEANLEF[z_]:=
  With[{f=Plus@@W E^(logLMA1F[#1,z]&/@X)},
    GETMEANLEF[z,W,X,f]];					
\end{verbatim}

Our tests indicate acceptable accuracy and about a $70\%$ improvement in speed
as compared with {\it Mathematica}'s more sophisticated {\bf NIntegrate}
function.

\subsection{Simulation Results for $n=50$}

\begin{figure}
\centerline{\psfig{figure=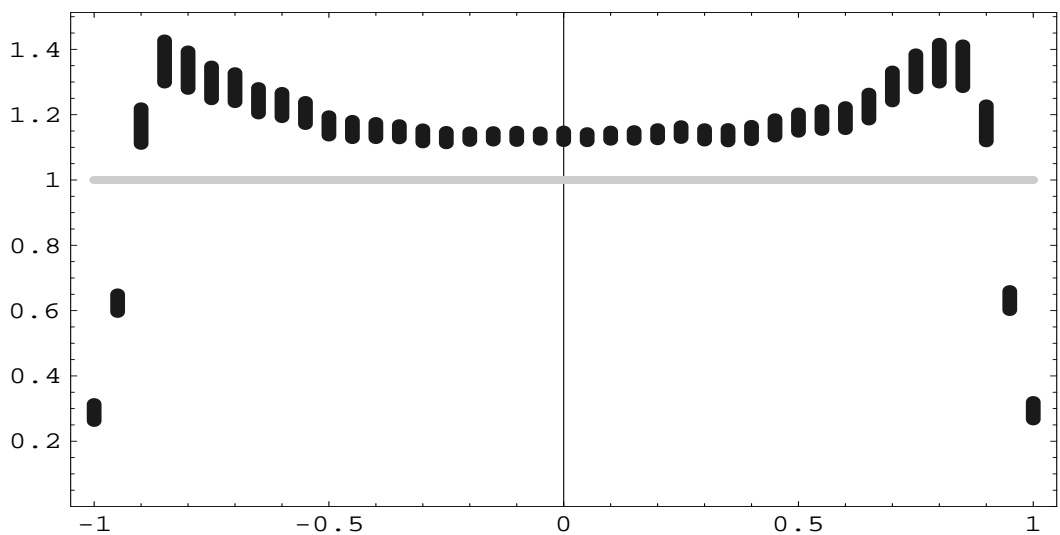,width=5.5in,height=2.75in}}
\caption{Empirical relative efficiency based on $10^4$ simulations of the MA(1) with $\mu=0$ and $n=50$\@. The length of the thick vertical lines indicate a $99.9\%$ confidence interval for $R(\bar{\theta}, \hat{\theta})$}
\label{Rma1n50mele}
\end{figure}

\begin{figure}
\centerline{\psfig{figure=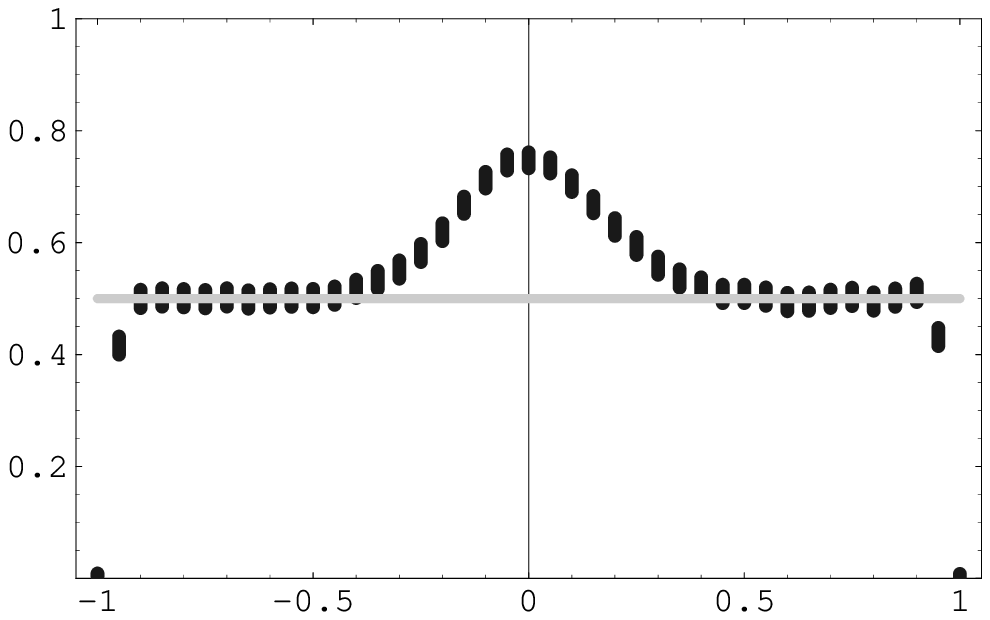,width=5.5in,height=2.75in}}
\caption{Empirical Pitman measure of closeness based on $10^4$ simulations of the MA(1) with $\mu=0$ and $n=50$\@. The length of the thick vertical lines indicate a $99.9\%$ confidence interval for ${\rm PMC}(\bar{\theta}, \hat{\theta})$}
\label{PMCma1n50mele}
\end{figure}

Using the {\it Mathematica} algorithms for the \mle\ and \mele\ derived above,
we determined $99.9\%$ confidence intervals for $R(\bar{\theta}, \hat{\theta})$
and ${\rm PMC}(\bar{\theta}, \hat{\theta})$ based on $10^4$ simulations for
each of the $41$ parameter values
$\theta= -1, -0.95, -0.90, \ldots, 0.95, 1$\@.
Figures~\ref{Rma1n50mele} and \ref{PMCma1n50mele}
show that the \mele\ dominates except for the cases
$\theta=\pm 1, \pm 0.95$\@.
We can safely conclude that the \mele\ is a better overall estimator than the
\mle.
Of course, as already pointed out another cogent reason for
preferring the \mele\ to the \mle\ is that it does
not produce noninvertible models.

If prior information is available then even better results
can be obtained.
Marriott and Newbold (1998)
have developed an ingenious approach to the unit root problem
in time series by noting this fact.

The simulations were repeated with the mean $\mu$ estimated by the sample
average and there was no major change is results.
The reader may like compare the estimators for other values of $n$ using the
{\it Mathematica} functions available in the electronic supplement.

\begin{figure}
\centerline{\psfig{figure=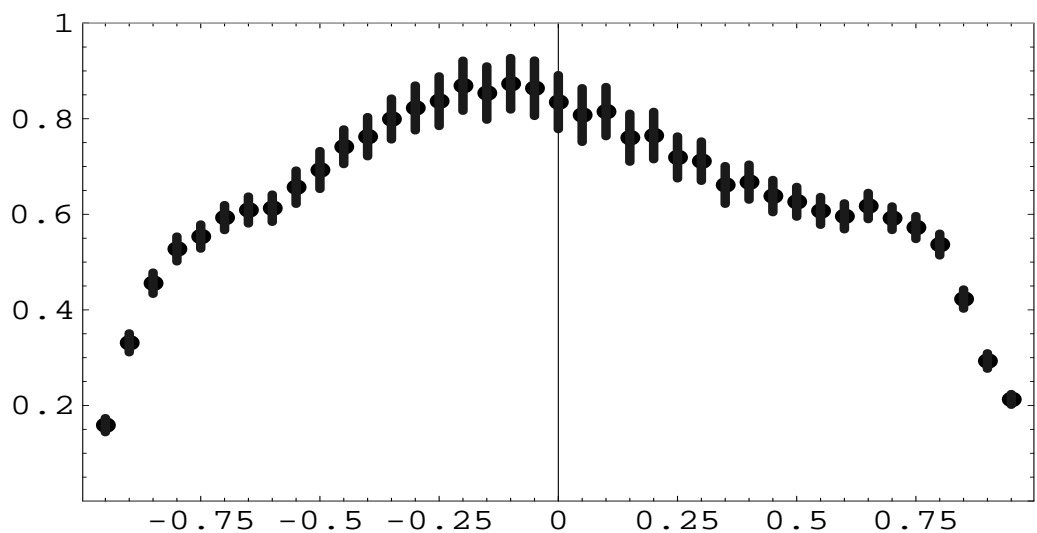,width=5.5in,height=2.75in}}
\caption{Empirical relative efficiency of Bayes estimate using a Jeffrey's prior.
Based on $10^4$ simulations of the MA(1) with $\mu=0$ and $n=50$\@.
The length of the thick vertical lines indicate a $99.9\%$ confidence
interval for $R(\tilde{\theta}, \hat{\theta})$}
\label{Rma1n50bayes}
\end{figure}

\begin{figure}
\centerline{\psfig{figure=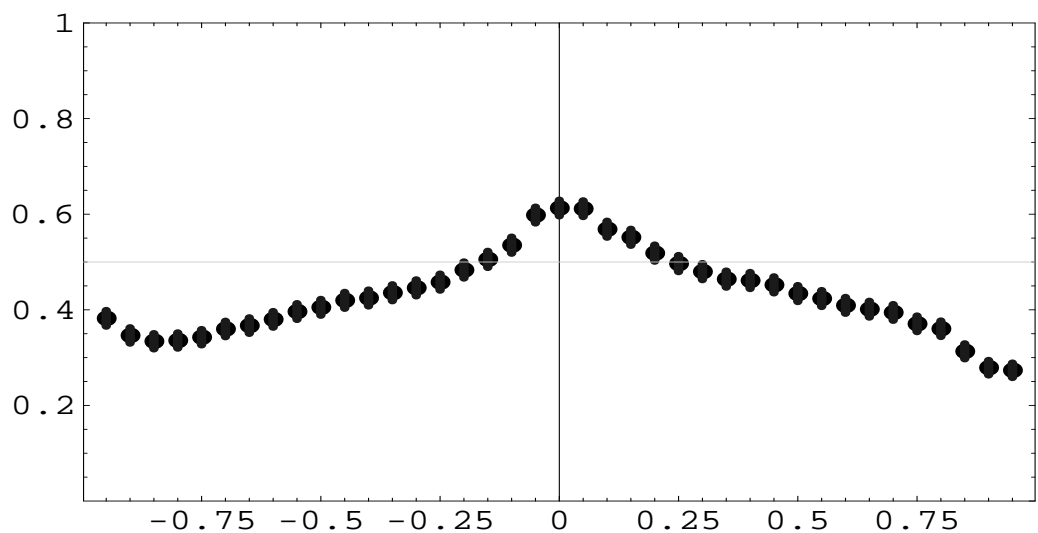,width=5.5in,height=2.75in}}
\caption{Empirical Pitman measure of closeness of Bayes estimate using a Jeffrey's prior.
Based on $10^4$ simulations of the MA(1) with $\mu=0$ and $n=50$\@.
The length of the thick vertical lines indicate a $99.9\%$ confidence interval
for ${\rm PMC}(\tilde{\theta}, \hat{\theta})$}
\label{PMCma1n50bayes}
\end{figure}

In the standard Bayesian analysis of the MA(1) model
the prior is given by
(Box and Jenkins, 1976, p. 250--258)
$$\pi(\theta) = 1/\sqrt{1-\theta^2}.$$
The computations were repeated using this prior and as shown in
Figures~\ref{Rma1n50bayes} and \ref{PMCma1n50bayes} the Bayes estimate with a
Jeffrey's prior performs about the same as the \mele.

\section{Concluding Remarks}

Previously Copas (1966) found that for AR(1) models, the \mele\ had lower \mse\
over much of the parameter region.
Our results show that for the MA(1) the improvement is even somewhat better.
The \mse\ is lower over a broader range and the piling-up effect on the \mle\ is
avoided.
Quenneville (1993) investigated the small sample properties of the \mele\ for
many other time series models and gave a general algorithm for the \mele\ in ARMA
models and found that in many cases the \mele\ produced estimates with smaller
\mse\ over most of the parameter region.
This work is further extended to state space prediction in
Quenneville and Singh (1999).

We would also like to mention that in our opinion {\it Mathematica\/} provides
an excellent and indeed unparalled environment for many types of
fundamental mathematical statistical research.
In comparison, no other computing environment provides such high quality
capabilities {\it simultaneously\/} in:
symbolics, numerics, graphics, typesetting
and programming.
The importance of a powerful user-oriented programming language
for researchers is sometimes lacking in other environments.
Stephan Wolfram once said that in his opinion the APL computing
language had many good ideas in this direction and that {\it Mathematica\/} has
incorporated all these capabilities and much more.
A partial check on this is given in
McLeod (1999) where it was found that most APL idioms
could be more clearly expressed in {\it Mathematica\/}.

However, for applied statistics and data analysis, Splus may
still be advantageous due
to the wide usage by researchers and the high quality functions for
advanced statistical methods that are available with Splus and
in the associated infrastructure.
From the educational viewpoint though this advantage may not be
so important since many students and researchers like to understand
the principles involved and with {\it Mathematica\/}
it is as easy to write out the necessary functions in {\it Mathematica\/}
notation as it would be to explain the procedures in a traditional mathematical
notation.

\section{References}
\begin{enumerate}
\item Andrews, D.F. (1997) Asymptotic expansions of moments and cumulants.
{\it Technical Report}.
\item Andrews, D.F. and Stafford, J.E. (1993)
Tools for the symbolic computation of asymptotic expansions.
{\it Journal of the Royal Statistical Society, B 55}, 613--627.
\item Barnard, G.A., Jenkins, G.M. and  Winsten, C.B.  (1960)
Likelihood inference and time series
{\it Journal of the Royal Statistical Society Series\/}  A 125,  321--372.
\item Borowski, E.J. and Borwein, J.M. (1991)
{\it The HarperCollins Dictionary of Mathematics}. New York: HarperCollins.
\item Box, G.E.P. and Jenkins, G.M. (1976)
{\it Time Series Analysis: Forecasting and Control\/}.
San Francisco: Holden-Day.
\item Box, G.E.P. and Tiao, G.C. (1973)
{\it Bayesian Inference in Statistical Analysis\/}.
Reading: Addison-Wesley.
\item Brockwell, P.J. and Davis, R.A. (1991)
{\it Time Series Theory and Methods}. New York: Springer-Verlag.
\item Cabrera, J. F. (1989)
Some experiments with maximum likelihood estimation using symbolic
manipulations.
{\it Proceedings of the 21st Symposium on the Interface of Statistics and
Computer Science}, Edited by K. Berk.
\item Copas, J. B. (1966)  Monte Carlo results for estimation in a stable
Markov time series.
{\it Journal of the Royal Statistical Society, A 129}, 110--116.
\item Cryer, J.D. and Ledolter, J. (1981)
Small-sample properties of the maximum likelihood estimator in the first-order
moving average model. {\it Biometrika, 68}, 691--694.
\item Currie, I.D. (1995) Maximum likelihood estimation and Mathematica.
{it Applied Statistics, 44}, 379--394.
\item Dempster, A.P. (1998)
Logicist statistics I. Models and modeling.
{\it Statistical Science\/} 13, 248--276.
\item Johnson, N.I. and Kotz, S. (1972)
{it Distributions in Statistics: Discrete Distributions}. New York: Wiley.
\item Keating, J.P., Mason, R.L. and Pranab, K.S. (1993)
{\it Pitman's Measure of Closeness}. Philadelphia: SIAM.
\item Newbold, P. (1974)
The exact likelihood function for a mixed autoregressive-moving average
process. {\it Biometrika, 61}, 423--426.
\item McLeod, A.I. (1999)
{\it APL\/}-idioms in {\it Mathematica\/}.
http://www.stats.uwo.ca/mcleod/epubs/idioms.
\item McLeod, A.I. and Quenneville, B. (1999)
{\it Mathematica\/} notebooks to accompany
Mean likelihood estimation.
http://www.stats.uwo.mcleod/epubs/mele.
\item Marriott, J. and Newbold, P. (1998)
Bayesian comparison of ARIMA and stationary ARMA models,
{\it International Statistical Review\/} 66, 323--336.
\item Pitman, E.J.G. (1937)
The closest estimates of statistical parameters,
{\it Proceedings of the Cambridge Philosophical Society\/} 33, 212--222.
{\it Biometrika\/} 30, 391--421.
\item Pitman, E.J.G. (1938)
The estimation of the location and scale parameters of a continuous
population of any given form,
{\it Biometrika\/} 30, 391--421.
\item Quenneville, B. (1993) {\it Mean Likelihood Estimation and Time Series
Analysis}. Ph.D. Thesis, University of Western Ontario.
\item Quenneville, B. and Singh, A.C. (1999)
Bayesian prediction MSE for state space models with estimated parameters.
{\it Journal of Time Series Analysis\/} (to appear).
\item Ritov, Y. (1990)
Asymptotic efficient estimators of the change point with unknown distributions.
{\it The Annals of Statistics\/} 18, 1829--1839.
\item Rubin, H. and Song, K.S. (1995)
Exact computation of the asymptotic efficiency of the maximum likelihood
estimators of a discontinuous signal in a Gaussian white noise.
{\it The Annals of Statistics\/} 23, 732--739.
\item Stafford, J.E. and Andrews, D.F. (1993)
A symbolic algorithm for studying adjustments to the profile likelihood.
{\it Biometrika, 80}, 715--730.
\item Stafford, J.E., Andrews, D.F. and Wang, Y. (1994)
Symbolic computation: a unified approach to studying likelihood.
{\it Statistics and Computing, 4}, 235--245.
\item Tanner, M.A. (1993)
{\it Tools for Statistical Inference}. New York: Springer-Verlag.
\item Wolfram, S. (1996) {\it Mathematica}. Champaign: Wolfram Research Inc.
\end{enumerate}

\newpage
\thispagestyle{empty}

\end{document}